\documentstyle[amssymb,amsfonts,12pt]{article}
\newtheorem{theorem}{\sc Theorem}

\newtheorem{lemma}[theorem]{\sc Lemma}
\newtheorem{corollary}[theorem]{\sc Corollary}

\newcommand{\proof}{\noindent {\sc Proof. }} 
 
\def\rest{\mathord{\restriction}}
\newcommand{\force}{\Vdash}
\newcommand{\open}{\Bbb}

\def\deq{\mathop=\limits^{\rm def}} 
 
\newlength{\labparwidth}
\newcommand{\rge}{\mbox{\rm rge}}

\newcommand{\se}{\subseteq}
\newcommand{\set}[2]{\{#1 \colon #2\}} 
\newcommand{\fin}{$\Box$\par\medskip} 
\newcommand{\dmd}{\diamondsuit}

\newcommand{\ga}{\alpha}
\newcommand{\gb}{\beta}

\newcommand{\gd}{\delta}

\newcommand{\gk}{\kappa}

\newcommand{\go}{\omega}
 
\newcommand{\ha}{\aleph}
 
\newcommand{\oP}{{\open P}}
\newcommand{\oQ}{{\open Q}}

\newcommand{\has}{\ha_1\mbox{-separable}}
\newcommand{\cha}{cardinality $\ha_1$ }

\begin{document}
\vspace*{1in}
\centerline{{\bf ON COHERENT SYSTEMS OF PROJECTIONS}}
\centerline{{\bf FOR $\ha_1$-SEPARABLE GROUPS}}
\vspace{.25 in}
\centerline{Paul C. Eklof}
\smallskip
\centerline{University of California, Irvine}
\centerline{Irvine, CA, 92717, USA}
\vspace{.25in}
\centerline{Alan H.
Mekler}
\smallskip
\centerline{Simon Fraser University}
\centerline{Burnaby, B.C. V5A 1S6, Canada}
\vspace{.25in}
\centerline{ Saharon Shelah}
\smallskip
\centerline{Hebrew University}
\centerline{Jerusalem, Israel
}

\vspace{.5in}

\noindent
{\small 
{\bf Abstract.} It is proved consistent with
either CH or $\neg$CH that there is an $\ha_1$-separable group of
cardinality $\ha_1$ which does not have a coherent system of
projections. It had previously been shown that it is consistent with $\neg$CH
that every $\ha_1$-separable group of cardinality $\ha_1$
does have a
coherent system of projections.}

\bigskip

\section{Introduction}

An abelian group $A$ is called $\has$ if every countable subset of
$A$ is contained in a countable free direct summand of $A$.
 An $\has$ group which is not free was first constructed by Griffith
\cite{G}, extending a construction by Hill \cite{H} for torsion groups.
 Such groups have been extensively studied, for example, in
\cite{M1}, \cite{E}, \cite{M2}  and \cite{EM}. To show that a group $A$
is $\has$ it suffices to produce an unbounded set of projections onto
countable free subgroups, that is, a family $\{\pi_i \colon i
\in I\}$ of functions $\pi_i \colon A \rightarrow H_i$ such that
$\pi_i \circ \pi_i = \pi_i, H_i = \rge(\pi_i)$ is a countable free
group, and such that for every countable subset $X$ of $A$, there is
$i \in I$ with $X \subseteq H_i$. (In fact, the existence of such a family 
is obviously equivalent to saying that $A$ is $\has$.)

In most cases, 
the construction of 
 an $\has$ group $A$  yields a group with
 a stronger property: it has a
{\it coherent unbounded system of projections}, i.e., a family $\{\pi_i \colon i
\in I\}$ as above with the additional property that if $H_j \subseteq
H_i$, then $\pi_j \circ \pi_i = \pi_j$. In fact, one cannot prove in
ZFC that an $\has$ group of \cha fails to have this
stronger property, because Mekler \cite{M2} has shown that PFA + $\neg$CH
implies that every $\has$ group  of \cha has this property (and more: it is
 in standard form). 

It has also been shown that the question of whether an $\has$ group
has a coherent system of projections (in an apparently stronger sense
--- ``with respect to a filtration'' ---
to be defined below), is relevant to the study of dual groups.
Specifically, every $\has$ group, $A$, of \cha  which has a
coherent system of projections with respect to a filtration and is
such that $\Gamma(A) \neq 1$ is a dual group. 
(See \cite[XIV.3.1]{EM}. It is an open question whether it is provable in
ZFC that every $\has$ group of \cha is a dual group.)

Thus it is a natural question to ask whether or not it is
provable in ZFC that every $\has$
group (of cardinality $\ha_1$) has a coherent system of projections.
This is posed as an open question in \cite{EM}. Here we answer that
question in the negative by showing that it is consistent both with
CH and with $\neg$CH that there is an $\has$ group of cardinality
$\ha_1$ with no coherent unbounded system of projections. Moreover, such
a group can be constructed to have any
desired Gamma invariant (other than $0$) and to be filtration-equivalent to
an $\has$ group
which {\it does} have a coherent system of projections.

\section{Preliminaries}

We will generally adhere to the terminology and notation of \cite{EM}. All groups referred to will be of cardinality at most 
$\aleph _1$. A {\it filtration} of an $\aleph _1$-separable group $A$ is 
 a continuous chain $\{A_\nu \colon \nu  < \omega _1\}$ of
subgroups of $A$ such that  $A_0 = 0$,  $A 
= \bigcup _{\nu <\omega _1}A_\nu $, and for all  $\nu  < \omega _1$, $A_{\nu +1}$ 
is a countable free direct summand of $A$. A homomorphism $\pi\colon  A
\rightarrow A$ is a {\it projection} if $\pi^2 = \pi$; in that case,
the image, $H$, of $\pi$ is a direct summand of $A$.

Given an $\aleph_1$-separable group $A$ and a filtration $\{A_\nu
\colon \nu \in \go_1\}$ of $A$, let $$
E \deq 
 \{\nu  \in  \lim (\omega _1)\colon A_{\nu +1}/A_\nu \hbox{ is not free}\}.
$$
Define $\Gamma(A) = \tilde{E}$, the equivalence class of $E$ modulo
the closed unbounded filter on ${\cal P}(\go_1 )$ (cf. \cite[II.4.4
and IV.1.6]{EM}).

A {\it coherent system of projections with respect to the filtration}
$\{A_\nu
\colon \nu \in \go_1\}$ of $A$
is a family of
projections $\{\pi _\nu { }\colon A \rightarrow  A_\nu  \colon  \nu
\notin  E \}$ 
such that for all  $\nu  < \tau $  in $\omega _1 \setminus  E$,  
$\pi _\nu { } \circ  \pi _\tau  = \pi _\nu { }.$ 

Clearly, $\{\pi _\nu { }\colon A \rightarrow  A_\nu  \colon  \nu
\notin  E \}$ is a coherent unbounded system of projections, as
defined in the Introduction. We do not know if, conversely, any
$\has$ group which has a coherent unbounded system of projections
also has a coherent system of projections with respect to a filtration.

We say that an $\has$ group $A$ has {\it quotient type H} if 
$A$ has a filtration $\{A_\nu
\colon \nu \in \go_1\}$ such that $A_{\nu + 1}/A_\nu \cong H$ for all
$\nu$ such that $A_{\nu + 1}/A_\nu$ is not free. (See \cite[p. 251]{EM}.)

Let $\mbox{\rm succ}(\go_1)$ (respectively, $\mbox{\rm lim}(\go_1)$)
denote the set of all successor (resp., limit) ordinals in $\go_1$.

\section{Construction of a counterexample using $\diamondsuit$}

	For a prime $p$, ${\open Q}^{(p)}$ denotes the subgroup of
$\open Q$ consisting of rationals whose denominators are a power of $p$.

\begin{theorem}
\label{diamond}
   Assume $\diamondsuit _{\omega _1}(S)$, where $S$ is a stationary set
of limit ordinals $< \omega _1$. Let $p$ be a prime. Then there exists an 
$\aleph _1$-separable group $A$ of cardinality $\aleph _1$ such that 
$\Gamma(A) = \tilde{S}$, $A$ is of 
quotient type ${\open Q}^{(p)}$, and $A$ has no coherent unbounded
system of projections. 
\end{theorem}

\proof   Let $D$ be the ${\open Q}$-vector space with basis $\{x_{\nu
,n}\colon  n \in  \omega $, $\nu  < \omega _1\} \cup  \{y_\delta
\colon  \delta  \in  S\}$. Let $D_\alpha$  be the subspace of $D$
generated by $\{x_{\nu ,n}\colon  n \in  \omega $, $\nu  < \alpha \}
\cup  \{y_\delta \colon  \delta  \in  S \cap  \alpha \}$. We shall
define inductively subgroups $A_\alpha $ of $D_\alpha $ such that for
all $\mu \geq  \alpha $, $A_\mu  \cap  D_\alpha  = A_\alpha $. At the
same time, we will define homomorphisms $t_{\alpha \nu }\colon
A_\alpha  \rightarrow  A_\nu $ for all successor ordinals $\nu  <
\alpha $. Our inductive construction will satisfy: 
\begin{quote}
(1) for all successor ordinals $\nu $ and all  $\gamma  > \alpha  > \nu $, 
$A_\nu$ is free and $t_{\alpha \nu }\rest A_\nu $ is the identity
(i.e., $t_{\alpha \nu }$ is a  
projection onto $A_\nu$) and $t_{\gamma \nu }\rest A_\alpha  =
t_{\alpha \nu };$  
\end{quote}
\begin{quote}
(2) if $\alpha  \notin  S$, then $A_{\alpha +1}/A_\alpha $ is free and if 
$\alpha  \in  S$, then $A_{\alpha +1}/A_\alpha  \cong  {\open Q}^{(p)}$. 
\medskip
\end{quote}
When the construction is completed we will define $A =\cup _{\alpha
<\omega _1} 
A_\alpha $ and 
$$t_\nu  = \cup _{\alpha < \omega _1}t_{\alpha \nu }\colon  A 
\rightarrow  A_\nu $$
 for each successor ordinal $\nu  < \omega _1$. We will 
carry out the construction so that the following properties will hold: 
\begin{quote}
(I) for every projection $\pi \colon A \rightarrow  H$ onto a
countable subgroup $H$ of $A$, there is a finite set $W_\pi  \subseteq
\mbox{\rm succ}(\omega _1)$ such that for all $a \in  A$, if $t_\nu
(a) = 0$ for all $\nu  \in  W_\pi $, then $\pi (a) = 0.$ 
\end{quote}
\begin{quote}
(II) whenever $W_0$ and $W_1$ are finite subsets of $\mbox{\rm
succ}(\omega _1)$ and $\beta  = \sup (W_0 \cap  W_1)$, there exists
$\delta  > \beta $ and $y_{\delta ,0}$, $y_{\delta ,1} \in  A_{\delta
+1}$ such that $0 \neq py_{\delta ,1} - y_{\delta ,0} \in  A_{\beta
+1}$, and $t_\nu (y_{\delta ,\ell}) = 0$ for all $\nu  \in  W_\ell$
($\ell = 0$, 
1).  
\end{quote}
Suppose for a moment that we can carry out the construction. Then $A$
is $\aleph _1$-separable since $\{ t_\nu \colon \nu \in
\mbox{succ}(\go_1) \}$ is an unbounded system of projections. Also,
(2) implies that $\Gamma(A) = \tilde{S}$ and $A$ has quotient type
${\open Q}^{(p)}$.

  We claim that there is no coherent unbounded
system of projections. Suppose, to the
contrary that $\{\pi_i \colon  i \in I\}$
is a coherent unbounded system of projections where $\rge(\pi_i) = H_i$.
 Then by (I), for each
$\pi _i $ there is a finite set $W_i
$ such that for all $a \in  A$, if $t_\nu (a) = 0$ for all $\nu  \in
W_i $, then $\pi _i (a) = 0$. Now apply the $\Delta $-system Lemma
\cite[p. 225]{J}: there is a finite set $\Delta  \subseteq  \omega _1$
and an uncountable subset $Z$ of $I$ such
that for all $i  \neq  i '$ in $Z$, $W_i  \cap W_{i '} =
\Delta $. Let $\beta  = \sup (\Delta )$. Choose $i_0$, $i_1 \in
Z$ such that $A_{\beta + 1}  \subseteq  H_{i _0}$ and $H_{i _0} \subseteq
 H_{i_1}$.
 Let $\delta $ and $y_{\delta ,0}$ and $y_{\delta ,1}$ be as in
(II) for $W_{i_0}$ and $W_{i_1}$. Then by (I) and (II) we have $\pi_{i
_\ell}(y_{\delta ,\ell}) = 0$ for $\ell = 0, 1$. By coherence we then
have $\pi_{i _0}(y_{\delta ,1}) = \pi _{i_0}(\pi _{i _1}(y_{\delta ,1}))
= 0$, so $\pi _{i _0}(py_{\delta ,1} - y_{\delta ,0}) = 0$, which is
a contradiction because $py_{\delta ,1} - y_{\delta ,0}$ is non-zero
and belongs to  $A_{\beta
+ 1}  \subseteq  H_{i _0}.$

\medskip

So it remains to do the construction. First let us write $S$ as the disjoint 
union 
$$
S = S_0 \amalg  S_1
$$

\noindent
of (stationary) sets such that $\diamondsuit _{\omega _1}(S_i)$ holds
for $i = 0$, 1. Also, choose a surjection $\psi $  from $S_0$ onto the
set of all pairs $(W_0$, $W_1)$ of finite subsets of $\mbox{\rm
succ}(\omega _1)$ such that for each $\delta  \in  S_0$, if $\psi
(\delta ) = (W_0$, $W_1)$, then $\delta  > \sup (W_0 \cap  W_1) +
\omega .$ 

Suppose now that we have constructed $A_\alpha $ and $t_{\alpha \nu }$
for all  $\alpha  < \gamma $. There are four cases to consider. 

In the first case, $\gamma $ is a limit ordinal. In this case, we let
$A_\gamma = \cup _{\alpha <\gamma } A_\alpha $ and $t_{\gamma \nu } =
\cup _{\nu <\alpha <\gamma } t_{\alpha \nu }$ for all successor
ordinals $\nu < \gamma $. Clearly (1) and (2) are satisfied. So now we
can assume that $\gamma  = \delta  + 1$ for some $\delta .$ 

In the second case, $\delta  \notin  S$. In this case we let $A_\gamma
= A_\delta  \oplus  \bigoplus _{n\in \omega } {\open Z}x_{\delta ,n}$
and for each successor $\nu  \leq  \delta $ we define $t_{\gamma \nu
}$ to be an extension of $t_{\delta \nu }\colon  A_\delta  \rightarrow
A_\nu $ (where $t_{\delta \delta }$ is the identity map if $\delta
\notin \mbox{\rm
succ}(\omega _1)$)
 such that the
$t_{\gamma \nu }$ $(\nu  \in
\mbox{\rm succ}(\omega _1) \cap  \gamma )$ satisfy: 
\begin{quote}
(3) $t_{\gamma \nu }(x_{\delta ,0}) = 0$ and for every finite subset
$F$ of $\mbox{\rm succ}(\omega _1) \cap 
\gamma $ and function $\theta \colon  F \rightarrow  A_\delta $, there
exists $k \geq  1$ such that $t_{\gamma \nu }(x_{\delta ,k}) = \theta
(\nu )$ for all $\nu  \in  F$, and $t_{\gamma \nu }(x_{\delta ,k}) =
0$ for $\nu  \notin  F.$
\end{quote}
\noindent
Since the number of pairs $(F$, $\theta )$ is countable, this is easy to 
arrange. 

In the third case, $\delta  \in  S_0$. Here we will do the construction to 
insure that (II) holds. Let $\psi (\delta ) = (W_0$, $W_1)$ and let $\beta  = 
\sup (W_0 \cap  W_1)$. Choose a ladder $\eta $ on $\delta $ such that
$\eta (0) = \beta $ and $\eta (n)$ is a successor ordinal greater than
$\sup (W_0 \cup  W_1)$ for all $n \geq  1$. By (3) there exists $k_1$
such that 
$$
t_{\delta \nu }(x_{\eta (1),k_1}) = -t_{\delta \nu }(x_{\eta (0),0})
$$

\noindent
for all $\nu  \in  W_0 \setminus  W_1$ and
$$
t_{\delta \nu }(x_{\eta (1),k_1}) = 0
$$

\noindent
for all other successor $\nu  \leq  \eta (1)$ (hence for all $\nu  \in
W_1)$.  

Now let $a_0 = px_{\eta (0),0}$, $a_1 = x_{\eta (1),k_1}$ and $a_j = 
x_{\eta (j),0}$ for $j \geq  2$. Let 
$$
y_{\delta ,n} = (y_\delta  + \Sigma _{j<n} p^ja_j)/p^n \in  D_{\delta +1}
$$

\noindent
(so $y_{\gd,0} = y_\gd$). Let $A_{\delta +1} = A_\gamma $ be the
subgroup of $D_{\delta +1}$ generated by
$$
A_\delta  \cup  \{y_{\delta ,n}\colon  n \in  \omega \}.
$$

\noindent
For all successor $\nu < \gd$ let
$$
t_{\gamma \nu }(y_{\delta ,n}) = - \Sigma ^\infty_{j=n}
p^{j-n}t_{\delta \nu }(a_j)
$$
for all $n \in \go$. This is easily seen to be a finite sum, by our
choice of the $a_j$, and the projections are well-defined. Moreover,
for $\nu  \in  W_1 \setminus  W_0$,

$$
t_{\gamma \nu }(y_{\delta,0}) = -p\cdot t_{\delta \nu }(x_{\eta (0),0})
$$

\noindent
and $t_{\gamma \nu }(y_{\delta ,n}) = 0$ for $n \geq  1$. For $\nu  \in  W_0 
\setminus  W_1$,
$$
t_{\gamma \nu }(y_{\delta ,1}) = t_{\delta \nu }(x_{\eta (0),0})
$$

\noindent
and $t_{\gamma \nu }(y_{\delta ,n}) = 0$ for $n \neq  1$. For $\nu  \in  W_0 
\cap  W_1$, since  $\nu  \leq  \beta  = \eta (0)$,
$t_{\gamma\nu}(x_{\eta(n),0}) = 0$ by definition; hence
$t_{\gamma \nu }(y_{\delta ,n}) = 0$ for all  $n$. Note also that 
$$py_{\gd,1} - y_{\gd,0} = a_0 = px_{\eta(0),0} = px_{\gb,0} \in
A_{\gb +1}.$$
\noindent
Hence, (II) is satisfied.

In the fourth and last case, $\delta  \in  S_1$. Then $\diamondsuit
(S_1 )$ gives us a prediction of a function $\pi _\delta \colon
A_\delta  \rightarrow  A_\delta $. If $\pi _\delta $ is not a
projection, or if there is a finite subset $W$ of $\mbox{\rm
succ}(\omega _1) \cap  \delta $ such that for all $a \in  A_\delta $,
$t_{\delta \nu }(a) = 0$ for all $\nu  \in  W$ implies $\pi _\delta
(a) = 0$, then define $A_\gamma $ and $t_{\gamma \nu }$ in any way
that satisfies (1) and (2). Otherwise, we want to define $A_\gamma $
so that, in addition, $\pi _\delta $ does not extend to $A_\gamma $.
Now $\pi _\delta $ is a projection:  $A_\delta  \rightarrow  H$ (for
some countable $H = \rge(\pi _\delta ))$ and if we write $\mbox{\rm
succ}(\omega _1) \cap  \delta $ as the increasing union, $\cup _{n \in
\omega } W_n$, of finite sets, then for each $n \in  \omega $ there
exists $a_n \in  A_\delta $ such that $\pi _\delta (a_n) \neq  0$ but
$t_{\delta \nu }(a_n) = 0$ for all $\nu  \in  W_n$. By the Lemma
following, there is a choice of $c_n \in  {\open Z}$ such that the
sequence $\langle \Sigma ^n_{j=0} p^jc_j\pi _\delta (a_j) \colon  n
\in  \omega \rangle $ does not have a limit in $H$ (in the $p$-adic
topology). Define
$$
y_{\delta ,n} = (y_\delta  + \Sigma _{j<n} p^jc_ja_j)/p^n \in  D_{\delta +1}
$$

\noindent
and let $A_{\delta +1} = A_\gamma $ be the subgroup of $D_{\delta +1}$
generated by $A_\delta  \cup  \{y_{\delta ,n}\colon  n \in  \omega
\}.$ Define $$ t_{\gamma\nu}(y_{\delta,n}) = - \Sigma_{j \geq n} p^{j
-n} t_{\delta\nu}(a_j) $$ which is well defined since almost all the
$t_{\delta\nu}(a_j)$ are $0$.  Then $\pi _\delta $ does not extend to
a homomorphism $h\colon  A_\gamma  \rightarrow  H$ since if it did,
$h(y_\delta )$ would be a limit of $\langle \Sigma ^n_{j=0} p^jc_j\pi
_\delta (a_j) \colon  n \in  \omega \rangle .$ 

This completes the inductive construction. It remains to check that
(I) holds.  Given any projection $\pi \colon A \rightarrow  H$, by
the diamond property, there
is a stationary subset $S'$ of $S_1$ such that for $\delta  \in  S'$,
$\pi \rest A_\delta  = \pi _\delta $. Hence, since $\pi _\delta $ does
extend to $A_{\delta +1}$, there is a finite subset $W_\delta $ of
$\mbox{\rm succ}(\omega _1) \cap  \delta $ such that for all $a \in
A_\delta $, $t_{\delta \nu }(a) = 0$ for all $\nu  \in W_\delta $
implies $\pi (a) = 0$. Then by Fodor's Lemma (cf. \cite[II.4.11]{EM}) and a
coding argument, there is a finite set $W_\pi $ such that for a
stationary subset $S''$ of $S'$, $\delta  \in  S''$ implies $W_\delta
= W_\pi $. Since $S''$ is unbounded in $\omega_1$, we are done.   \fin

\begin{lemma}
\label{hat}
   Let $H$ be a countable free group and $\hat{H}$ its closure in the 
$p$-adic topology. If $\langle b_n\colon  n \in  \omega \rangle $ is
a sequence of non-zero elements of $H$, then 
$$
\{\Sigma _{j \in  \omega } p^jc_jb_j\colon  \langle c_j\colon  j \in
\omega \rangle  \in {\open Z}^\omega  \}
$$
\noindent
is a subset of $\hat H$ of cardinality $2^{\aleph _0}.$
\end{lemma}

\proof By induction choose an increasing sequence $(m_n)$, so that
$p^{m_n + n}$ does not divide any element of $\{p^{m_k + k}b_k\colon k <
n\}$. For any $\xi \in {}^\go2$ let $c_{\xi n} = \xi(n)p^{m_n}$. It
remains to check that if $\xi_0 \neq \xi_1$
then $\sum^\infty_{k = 0} 
p^kc_{\xi_0 k}b_k \neq \sum^\infty_{k = 0} p^kc_{\xi_1 k}b_k$. Let $n$
be minimal so that $\xi_0(n) \neq \xi_1(n)$, then 
$$\sum^n_{k = 0} p^kc_{\xi_0 k}b_k - \sum^n_{k = 0} p^kc_{\xi_1 k}b_k
= \pm p^{m_n + n} b_n \not\equiv 0 \bmod p^{m_{n+1} + n+1} H.$$
However, $p^{m_{n+1}+ n+1}$ divides $\sum^\infty_{k = n+1}
p^kc_{\xi_0 k}b_k - \sum^\infty_{k = n+1} p^kc_{\xi_1 k}b_k$. \fin

\begin{corollary}
It is consistent with ZFC that there are filtration-equivalent
$\ha_1$-separable groups $A$ and $B$ such that $B$ has a coherent
system of projections with respect to a filtration but $A$ does not
have a coherent unbounded system of projections.
\end{corollary}

\proof
Let $A$ be as constructed in the Theorem. Associated with each
$\delta \in S$ there is a ladder $\eta_\gd$ on $\gd$ such that $p^{n +1}$
divides $y_{\gd, 0}$ mod $A_\nu$ if and only if 
$\nu \geq \eta_\gd(n)$. If we construct $B$ as in \cite[VIII.1.1]{EM}
(with $p_\gd = p$ for all $\gd \in S$), then by \cite[VII.1.10]{EM} $B$
has a coherent system of projections with respect to a filtration
and by \cite[Thm. 1.4]{E}, $A$ 
and $B$ are filtration-equivalent.  \fin

The following should be compared with \cite[XIV.3.1]{EM}. (See also the
introductory remarks concerning dual groups.)

\begin{corollary}
It is consistent with ZFC that there is an $\ha_1$-separable group
$A$ such that $\Gamma(A) \neq 1$ and $A$ does not have a coherent
system of complementary summands.  \fin

\end{corollary}

\section{Counterexamples  where CH fails}

 Theorem~\ref{diamond} requires $\dmd(S)$ which implies CH. We know
that it is consistent with $\neg$CH that every $\has$ group of \cha
has a coherent unbounded system of projections (cf. \cite{M2}).  So
the question naturally arises whether it is consistent with $\neg$CH
that there is an $\ha_1$-separable group of cardinality $\ha_1$ which
does not have a coherent unbounded system of projections. Here we
shall prove that the answer to the question is ``yes". In fact the
forcing used is just the simplest possible, namely Fn($\gk, 2, \go$),
the forcing for adding $\gk$ Cohen reals, where $\gk \geq \ha_2$, to
make CH fail. (Fn($\gk, 2, \go$) is the poset consisting of all
partial functions from $\gk$ to $2$ whose domains have cardinality
less than $\go$.)

\begin{theorem}
    It is consistent with $\neg$CH that for every stationary
subset $S$ of $\lim(\go_1)$ there is an
$\ha_1$-separable group $A$ of cardinality $\ha_1$ with
$\Gamma(A) = \tilde{S}$
which does not have a
coherent unbounded system of projections.
\end{theorem}

\proof 
We shall prove the following lemma.

\begin{lemma}
\label{cohen}
Suppose $\oP = \mbox{Fn}(\ha_1, 2, \go)$ and suppose $S \in V$ is a
stationary subset of $\lim(\go_1)$. If $G$ is generic for $\oP$, then
in $V[G]$ there is an $\ha_1$-separable group $A$ of cardinality
$\ha_1$ with $\Gamma(A) = \tilde{S}$ which does not have a coherent
unbounded system of projections. 
\end{lemma}

Assume for the moment that the lemma is correct.  Let $\oP'$ be
Fn($\gk, 2, \go$) where $\gk \geq \ha_2$, and let $G'$ be generic for
$\oP'$.  Given any stationary set,
$S$, in the generic extension, $V[G']$, we recast the forcing as a
two-step iteration, say $\oP_0\times\oP_1$ with generic set $G' = G_0
\times G_1$, where $\oP_0$ adds some number of Cohen reals, $\oP_1$
adds $\ha_1$ Cohen reals and $S \in V[G_0]$. By Lemma~\ref{cohen}
there is an $\ha_1$-separable group $A$ of cardinality $\ha_1$ in
$V[G_0][G_1] = V[G']$ with $\Gamma(A) = \tilde{S}$ and with no
coherent unbounded system of projections.

Thus it remains to prove Lemma~\ref{cohen}.  We will describe an
iterated forcing which forces the existence of the desired $A$.  The
forcing will be an iteration of length $\go_1$. Afterward we will note
that the forcing is equivalent to adding $\ha_1$ Cohen reals.  We
follow the usual notation where at each step $\ga$ the iterate is
$Q_\ga$ and the result of the iteration up to $\ga$ is $\oP_\ga$. We
will let $G_\ga$ denote an arbitrary $\oP_\ga$-generic set and talk of
members of $V[G_\ga]$ where more correctly we should talk of 
$\oP_\ga$-names.

	As well as constructing the sequence $Q_\ga$ we will define a
sequence of groups $A_\ga$ and projections $t_{\ga\nu}$ where $A_\ga,
t_{\ga\nu} \in V[G_\ga]$ and the $A_\ga$'s, $t_{\ga\nu}$ are as in
Theorem~\ref{diamond} (except for properties (I) and (II) which we
will have to verify). The groups $A_\ga$ will be constructed to be
subgroups of $D_\ga \subseteq D$, as in Theorem~\ref{diamond}.
 By coding we can assume that the set underlying
$D$ is $\go_1$. As in the proof of Theorem~\ref{diamond}, partition
$S$ into two disjoint stationary subsets $S_0$ and $S_1$.

	The construction goes by cases. If $\ga \notin S_1$ then define
$Q_\ga$ to be trivial (the one element poset). The construction of
$A_{\ga+1}$ and $\set{t_{\ga+1\nu}}{\nu \leq \ga, \nu \in
\mbox{succ}(\go_1)}$ is as in Theorem~\ref{diamond} (i.e., as in the
second or third case). Of course the
construction of $A_\gd$ and $t_{\gd, \nu}$ is determined when $\gd$ is
a limit ordinal. 

	Suppose now that $\gd \in S_1$. We will work in $V[G_\gd]$ and
define $\oQ_\gd$. Then $Q_\gd$ will be the obvious $\oP_\gd$-name.
List as $(\ga_n\colon n < \go)$ the ordinals in $\gd \cap
\mbox{succ}(\go_1)$. The 
forcing $\oQ_\gd$ is defined to be the set of sequences of the form $(c_0, a_0,
\ldots, c_{n-1}, a_{n-1})$ where for all $m < n$, $c_m \in \{0, 1\}$,
$a_m \in A_\gd$ and if $j < m$ then $t_{\gd, \ga_j}(a_m) = 0$. 
$\oQ_\gd$ is ordered
by extension. A generic set for $\oQ_\gd$ can be identified with a
sequence of length $\go$.  Given a generic set $G_{\gd+1}$ for
$\oP_{\gd+1}$ and so a generic sequence $(c_j, a_j\colon j <\go)$ for
$\oQ_\gd$, let $$y_{\gd,n} = (y_\delta + \Sigma _{m<n} p^mc_ma_m)/p^n
\in D_{\delta +1} $$ Let $A_{\delta +1} = A_\gamma $ be the subgroup
of $D_{\delta +1}$ generated by $ A_\delta \cup \{y_{\delta ,n}\colon
n \in \omega \} $. The definition of the projections is as in
Theorem~\ref{diamond}; they are well-defined because for all $j \in
\go$, for all $m > j$,
$t_{\gd, \ga_j}(a_m) = 0$.

	In $V[G_{\go_1}]$, we let $A = \bigcup_{\ga < \go_1} A$ and
for every successor ordinal $\nu$, we let $t_\nu = \bigcup_{\gb > \nu}
t_{\gb\nu}$. We will observe that $\oP_{\go_1}$ is equivalent to
adding $\ha_1$ Cohen reals. In particular, the forcing is c.c.c.\ and
so $\go_1$ is preserved and $A$ is an $\has$ group of
\cha.
To see that $A$ is the desired group we have to check that property
(I) from Theorem~\ref{diamond} holds. (The construction guarantees
that property (II) holds for exactly the same reasons as in the proof
of Theorem~\ref{diamond}). The proof that $A$ satisfies property (I) is
contained in the following two lemmas.

\begin{lemma}
    Use the notation above. Suppose $\gd \in S_1$. Furthermore suppose
$\pi\in V[G_\gd]$ and $\pi$ is a projection from $A_\gd$ to $H$ so
that for every  finite set $w \se \set{\ga < \gd}{\ga \in
\mbox{succ}(\go_1)}$ there 
is  $a \in A_\gd$ such that  $t_{\gd\ga}(a) = 0$ for all $\ga
\in w$ and $\pi(a) \neq 0$. Then $\pi$ does not extend to a projection
from $A_{\gd+1}$ to $H$.
\end{lemma}

\proof  We will work in  $V[G_\gd]$. Fix some such $\pi$. It suffices to
show for all $a \in A_\gd$ that 
$$D_a \deq \set{q \in \oQ_\gd}{q \force \mbox{
``if } \hat{\pi} \mbox{ is an extension of } \pi \mbox{ to } A_{\gd+1}
\mbox{ then } \hat{\pi}({y_\gd})\neq a\mbox{''}}$$ 
is dense. Fix $a
\in A_\gd$ and consider any condition $(c_0, a_0, \ldots, c_{n-1},
a_{n-1})$. Choose $a_n$ so that $\pi(a_n)\neq 0$ and
$t_{\gd\ga_m}(a_n) = 0$ for all $m < n$. For some choice of $c_n \in
\{0, 1\}$, $\sum_{m=0}^n p^mc_m\pi(a_m) \neq a$. Since $A_\gd$ is
free, there is $k > n$ so that $\sum_{m=0}^n p^mc_m\pi(a_m) \not\equiv
a \bmod p^kA_\gd$. For $m$ so that $n < m < k$ let $c_m = 0$ and let
$a_m = 0$. Notice that if $b_i$ $(i \geq k)$ are any
elements of $A_\gd$ we have
$$\sum_{m=0}^{k-1} p^mc_m\pi(a_m) + \sum_{m=k}^\infty p^mb_m \equiv
\sum_{m=0}^n p^mc_m\pi(a_m) \not\equiv a \bmod p^kA_\gd.$$
Hence $(c_0, a_0, \ldots c_{k-1}, a_{k-1})$ belongs to $D_a$. \fin

(We could have replaced Lemma~\ref{hat} by an argument like that in the
preceding proof.)

\begin{lemma}
    Suppose $\pi\in V[G_{\go_1}]$ is a projection of $A$ to a
subgroup $H$. Then there is a closed unbounded set $C$ so
that for all $\ga \in C$, $\pi\rest A_\ga \in V[G_\ga]$. (We assume
here, as we have done tacitly above, that $G_\ga$ is the restriction of
$G_{\go_1}$ to $\oP_\ga$.)
\end{lemma}

\proof This is a standard fact for finite support iterations of
c.c.c.\ forcing, so we  will just 
sketch the argument. Take $\tilde{\pi}$ a name for $\pi$. For each
$\ga \in A$, take $X_\ga$ a maximal antichain of conditions so that for all
$q \in X_\ga$, there is $a_{q\ga}$ so that $q \force \tilde{\pi}(\ga)
= a_{q\ga}$. (Recall that the underlying set of $A$ is contained in
$\go_1$). Since $\oP$ is c.c.c., each $X_\ga$ is countable. Our cub
$C$ consists of $\set{\ga < \go_1}{\mbox{for all }\gb \in A_\ga, X_\gb
\se \oP_\ga \mbox{ and for all }q \in X_\gb, a_{q\gb} \in A_\ga}$. \fin 

	It remains to observe that $\oP_{\omega_1}$ is equivalent to
adding $\ha_1$ Cohen reals. The proof uses two pieces of folklore.
The first one that any countable poset with the property that any
element has two incompatible extensions is equivalent to the forcing
for adding a Cohen real. The second, which uses the first, is that an
iteration of length $\go_1$ such that each iterate is forced to be a
countable poset with the property that any element has two
incompatible extensions  is equivalent to adding $\ha_1$ Cohen reals.
A somewhat fuller explanation can be found in the proof of Lemma~1.5
of \cite{MS}. If we view $\oP_{\go_1}$ as the iteration of
$\{Q_\gd\colon \gd \in S_1\}$, then the second piece of folklore
applies. \fin

\section{Questions}

	One question that we do not know the answer to is whether or
not the existence of an $\ha_1$-separable group of cardinality
$\ha_1$ without a coherent unbounded system of  projections
 follows from CH
alone. (Presumably one would use weak diamond in such a proof.) To
put the question a different way,  
    is it consistent with CH that every $\ha_1$-separable group of
cardinality $\ha_1$ has a coherent unbounded system of
projections?

	Another question along the same lines is whether MA + $\neg$CH
implies that every $\ha_1$-separable group of cardinality $\ha_1$ has
a coherent unbounded system of projections. Since PFA implies MA + $\neg$CH, we
know that it is consistent with MA + $\neg$CH that every
$\ha_1$-separable group of cardinality $\ha_1$ has a coherent system
of projections with respect to a filtration. Our methods cannot be immediately translated over to
a model of MA + $\neg$CH, since we have built a group which is
filtration equivalent to a group with a coherent system of projections,
while under MA + $\neg$CH any two filtration equivalent $\has$ groups of
\cha\ are isomorphic (\cite{E})

Finally, there is the question of whether the existence of a coherent
unbounded system of projections for an $\has$ group $A$ of
\cha implies the existence of a coherent system of projections with
respect to a filtration of $A$. (It clearly implies the existence of
a filtration
$\{A_\nu
\colon \nu \in \go_1\}$ of $A$
and a coherent family of
projections $\{\pi _\nu { }\colon A \rightarrow  A_\nu  \colon  \nu
\in \mbox{succ}(\go_1) \}$; the problem is to define coherently
projections $\pi_\nu$ when $\nu$ is a limit ordinal not in $E$.)

\bigskip

\centerline{ACKNOWLEDGEMENTS}

\smallskip

The research for this paper was carried out, in part, while the
authors were visiting Rutgers University. We  thank Rutgers
University for its support. Mekler's research was partially
 supported by NSERC grant A9848.
Shelah's research was partially supported by the BSF
(United States-Israel Binational Science Foundation).  Publication 426.

\end{document}